\documentclass{article}
\usepackage[utf8]{inputenc}
\usepackage{amsfonts}
\usepackage{amssymb}
\usepackage{amsmath} 
\usepackage{amsthm}
\usepackage{hyperref}
\usepackage{hyphenat}
\usepackage{flushend}
\usepackage{geometry}
\theoremstyle{remark}

\title{Non-extendability of holomorphic functions with bounded or continuously extendable derivatives}
\author{D. Moschonas, V. Nestoridis}
\date{}
\begin{document}
\frenchspacing
\maketitle \begin{center} \textbf{Abstract}
\end{center}
We consider the spaces $H_{F}^{\infty}(\Omega)$ and $\mathcal{A}_{F}(\Omega)$ containing all holomorphic functions $f$ on an open set $\Omega \subseteq \mathbb{C}$, such that all derivatives $f^{(l)}$, $l\in F \subseteq \mathbb{N}_0=\{ 0,1,...\}$, are bounded on $\Omega$, or continuously extendable on $\overline{\Omega}$, respectively. We endow these spaces with their natural topologies and they become Fréchet spaces. We prove that the set $S$ of non-extendable functions in each of these spaces is either void, or dense and $G_\delta$. We give examples where $S=\varnothing$ or not. Furthermore, we examine cases where $F$ can be replaced by $\widetilde{F}=\{ l\in \mathbb{N}_0:\min F \leqslant l \leqslant \sup F\}$, or $\widetilde{F}_0= \{  l\in \mathbb{N}_0:0\leqslant l \leqslant \sup F\}$ and the corresponding spaces stay unchanged.\\\\
\textit{AMS classification number:} ???

\section{Introduction}
Suppose that $f$ is a function  defined on an interval $I$ and that its derivatives $f^{(a)}$ and $f^{(b)}$ are bounded, where $a$ and $b$ are natural numbers, $0\leqslant a<b<+\infty$. Let $l$ be a natural number, such that $a<l<b$ and consider the derivative $f^{(l)}$. A natural question is if $f^{(l)}$ is also bounded, or more generally what can be said about the growth of $f^{(l)}$ on $I$. This question has been investigated by several mathematicians, such as Landau, Kolmogorov, Hardy, Littlewood and others; see [3], [4], [5], [6]. In particular, if $I=\mathbb{R}$ or $I=(0,+\infty)$, then the boundedness of $f^{(a)}$ and $f^{(b)}$ imply the boundedness of $f^{(l)}$, $a<l<b$ ([2], [3]). It follows that, if $\Omega$ is an open subset of the complex plane $\mathbb{C}$ and if $f$ is a holomorphic function on $\Omega$, such that $f^{(a)}$ and $f^{(b)}$ are bounded on $\Omega$, then all the intermediate derivatives $f^{(l)}$, $a<l<b$, are also bounded on $\Omega$, provided that $\Omega$ is the union of open half-lines. For instance, $\Omega$ could be an open angle, or a strip, or the union of two meeting angles, or of two meeting strips, etc.\\\\
The above consideration leads us to consider the space:
\begin{align*}
H_{F}^{\infty}(\Omega)=\{f:\Omega\rightarrow \mathbb{C}\ \text{holomorphic}: f^{(l)} \ \text{is bounded on} \ \Omega,\ \text{for all} \ l \in F\},
\end{align*}
where $F$ is an arbitrary, non-empty subset of $\mathbb{N}_0=\{ 0,1,... \}$ and examine whether $H_{F}^{\infty}(\Omega)=H_{\widetilde{F}}^{\infty}(\Omega)$ or not, where $\widetilde{F}=\{ l\in \mathbb{N}_0:\min F\leqslant l \leqslant \sup F\}$. Indeed, if $\Omega$ is a union of open half-lines, then $H_{F}^{\infty}(\Omega)=H_{\widetilde{F}}^{\infty}(\Omega)$. We believe that this does not hold for the general open set, but we do not have a counter-example. Furthermore, we believe that a complete metric can be defined on the set of all open sets $\Omega$ (contained in a disc), so that for the generic open set $\Omega$, it holds $H_{F}^{\infty}(\Omega)\neq H_{\widetilde{F}}^{\infty}(\Omega)$. The space $H_{F}^{\infty}(\Omega)$, endowed with its natural topology, is a Fréchet space and thus Baire's Category theorem is at our disposal in order to prove some generic results.\\\\
In analogy to the space $H_{F}^{\infty}(\Omega)$, we consider the space:
\begin{align*}
\mathcal{A}_{F}(\Omega)=\{f:\Omega\rightarrow \mathbb{C} \ \text{holomorphic}: f^{(l)}\ \text{has a continuous extension on} \ \overline{\Omega},\ \text{for all} \ l \in F\},
\end{align*}
where the closure is taken in $\mathbb{C}$. This space, endowed with its natural topology, is also a Fréchet space and Baire's theorem can be applied in order to prove some generic results. Moreover, we use the completeness of these spaces and a result from [14] to prove that either every function is extendable, or generically every function is non-extendable. We give examples where each horn of the above dichotomy occurs.\\\\
We note that if $p\in \mathbb{N}_0=\{ 0,1,...\}$ and $F=\{0,1,...,p\}$, the spaces $H_{F}^{\infty}(\Omega)$ and $\mathcal{A}_{F}(\Omega)$ are denoted by $H_{p}^{\infty}(\Omega)$ and $\mathcal{A}_{p}(\Omega)$, respectively; these spaces have been studied extensively in [13] and elsewhere.\\\\
Finally, in the last section we present another dichotomy result regarding the space $H_{F}^{\infty}(\Omega)$, proven using a result from [15]. It states that if $F$ is a non-empty subset of $\mathbb{N}_0=\{ 0,1,...\}$, $\Omega$ is an open subset of $\mathbb{C}$ and $l\notin F$, then either for every function $f$ in $H_{F}^{\infty}(\Omega)$ the derivative $f^{(l)}$ is bounded, or generically for every function $f$ in $H_{F}^{\infty}(\Omega)$ the derivative $f^{(l)}$ is unbounded. Thus, either $H_{F}^{\infty}(\Omega)=H_{F\cup \{l\}}^{\infty}(\Omega)$, or $H_{F\cup \{ l \}}^{\infty}(\Omega)$ is meager in $H_{F}^{\infty}(\Omega)$. It remains open to find an example of an open set $\Omega\subseteq \mathbb{C}$ for which the equality $H_{F}^{\infty}(\Omega)=H_{\widetilde{F}}^{\infty}(\Omega)$ fails to hold, though we believe that such examples exist and that this phenomenon is valid for the generic open set $\Omega$. We also prove that for any unbounded open set $\Omega$ and any $F\subseteq \mathbb{N}_0$, generically for every function $f$ in $\mathcal{A}_{F}(\Omega)$, all derivatives $f^{(l)}$, $l\in \mathbb{N}_0$, are unbounded.

\section{Preliminaries}
We give the following definition:\\\\
\textbf{Definition 2.1:} Let $\Omega \subseteq \mathbb{C}$ be an open set and $f:\Omega\rightarrow \mathbb{C}$ be a holomorphic function. We say that $f$ is extendable (in the sense of Riemann surfaces) if there exist two open discs $D_1$ and $D_2$, such that $D_2\cap \Omega \neq\varnothing$, $D_2\cap \Omega^{c} \neq \varnothing$ and $D_1\subseteq \overline{D_1}\subseteq D_2\cap \Omega$ and a (bounded) holomorphic function $F:D_2\rightarrow \mathbb{C}$, such that $F|_{D_1}=f|_{D_1}$. Otherwise, we say that $f$ is non-extendable, or that it is holomorphic exactly on $\Omega$.\\\\
\textbf{Remark 2.2:} The reason for which $F$ can be chosen to be bounded is that, if needed, we can replace disc $D_2$ by another disc $D_2^{'}$ compactly contained in $D_2$, that is $\overline{D_2^{'}}\subseteq D_2$, such that $\overline{D_1} \subseteq D_2^{'}\cap \Omega$. Also, in Definition 2.1, one can replace $D_2$ with any non-empty domain $U \subseteq \mathbb{C}$ satisfying $U\cap \Omega \neq\varnothing$ and $U\cap \Omega^{c} \neq \varnothing$ and replace $D_2$ with a connected component $V$ of $U\cap \Omega$, resulting in an equivalent definition of extendability; see [14] for a proof of this fact.\\\\
In [14] we find the following theorem, proven using Baire's and Montel's theorems:\\\\
\textbf{Theorem 2.3:} Let $\Omega \subseteq \mathbb{C}$ be an open set and $H(\Omega)$ be the set of holomorphic functions on $\Omega$. Also, let $X(\Omega)\subseteq H(\Omega)$ be a topological vector space endowed with the usual operations $+, \cdot$, whose topology is induced by a complete metric. Suppose that the convergence $f_n \rightarrow f$ in $X(\Omega)$ implies the pointwise convergence $f_n(z) \rightarrow f(z)$, for all $z\in \Omega$. Then, there exists an $f \in X(\Omega)$ which is non-extendable, if and only if, for any two discs $D_1$ and $D_2$ as in Definition 2.1, there exists a function $f_{D_1,D_2}\in X(\Omega)$ so that the restriction $f_{D_1,D_2}|_{D_1}$ on $D_1$ does not possess a (bounded) holomorphic extension on $D_2$. If the previous assumptions hold, then the set $S=S_{X(\Omega)}=\{ f\in X(\Omega):f$ is non-extendable$\}$ is a dense and $G_\delta$ subset of $X(\Omega)$.\\\\
\textbf{Corollary 2.4:} The set $S_{X(\Omega)}$ from Theorem 2.3 is always a $G_\delta$ subset of $X(\Omega)$, because either $S_{X(\Omega)}=\varnothing$, or $S_{X(\Omega)}\neq \varnothing$ and then it is dense and $G_\delta$ from the previous theorem.\\\\
We also present the following geometric lemma which will be useful in section 5. The proof we provide is quite elementary:\\\\
\textbf{Lemma 2.5}: Let $\Omega \subseteq \mathbb{C}$ be an open set. If $\Omega$ is unbounded and convex, then it is the union of open half-lines.\\\\
\textit{Proof:} Fix a point $\alpha$ in $\Omega$. Since $\Omega$ is unbounded, there exists a sequence $(z_n)_n$ of points in $\Omega$, such that $0<|z_n-\alpha|\rightarrow +\infty$. Let $h_n=\frac{z_n-\alpha}{|z_n-\alpha|}$, then $|h_n|=1$, for all $n\in \mathbb{N}$. That is, every $h_n$ belongs to the compact set $S^1=\{z\in \mathbb{C}: |z|=1\}$, therefore we can extract a convergent subsequence of $(h_n)_n$. For simplicity and without any loss of generality, we can assume that $(h_n)_n$ itself converges to some $h \in S^1$, as $n \rightarrow \infty$. We will show that $\alpha+th\in \Omega$, for all $t>0$ (for $t=0$ it is trivial). Indeed, let $t>0$. Since $|z_n-\alpha|\rightarrow +\infty$, for large enough $n$ we have that $|z_n-\alpha|>2t$ and therefore $0<\frac{2t}{|z_n-\alpha|}<1$. Since $\alpha \in \Omega$, $z_n \in \Omega$ for all $n\in \mathbb{N}$ and $\Omega$ is convex, we deduce that:
\begin{align*}
w_n=\left ( 1-\frac{2t}{|z_n-\alpha|} \right )\alpha+\frac{2t}{|z_n-\alpha|} z_n \in \Omega
\end{align*}
 for large enough $n$. Observe that:
\begin{align*}
w_n=\alpha+2t\frac{z_n-\alpha}{|z_n-\alpha|}=\alpha+2th_n \rightarrow a+2th.
\end{align*}
Therefore, $w_n-2th \rightarrow \alpha$ and thus $2\alpha+2th-w_n \rightarrow \alpha$. Since $\alpha \in \Omega$ and $\Omega$ is an open set, we deduce that $2\alpha+2th-w_n \in \Omega$, for large enough $n$. So far, we have showed that $w_n \in \Omega$ and $2\alpha+2th-w_n \in \Omega$, for large enough $n$. For these points, the convexity of $\Omega$ implies the following:
\begin{align*}
\frac{1}{2}w_n+\frac{1}{2} \left (2\alpha+2th-w_n  \right )=\alpha+th\in \Omega.
\end{align*}
Thus, the closed half-line $\{ \alpha+th: t\geqslant 0 \}$ is contained in $\Omega$. In order to find an open half-line contained in $\Omega$, we pick an $r>0$ sufficiently small, so that the disc $D(\alpha,r)$ is contained in $\Omega$. Then, by extending the closed half-line $\{ \alpha+th: t\geqslant 0 \}$ towards point $\alpha$ by an open line segment of length $r$, we conclude that the open half-line $\{ \alpha+th: t>-r \}$, which is parallel to $h$, is contained in $\Omega$. This completes the proof.$\qed$\\\\
\\\textbf{Remark 2.6:}
The result obtained in the previous lemma can be strengthened in the follo-  wing sense: if $\Omega \subseteq \mathbb{C}$ is an open, unbounded, convex set, then it is the union of open half-lines which are parallel to each other. This implication follows directly from the previous proof. Indeed, pick any two points $\alpha$, $\beta \in \Omega$ and consider $u_n=\frac{z_n-\alpha}{|z_n-\alpha|}$ and $v_n=\frac{z_n-\beta}{|z_n-\beta|}$, for $n\in \mathbb{N}$, where $(z_n)_n$ is a sequence of points in $\Omega$, such that $z_n \rightarrow \infty$ and $z_n \neq \alpha,\beta$ for all $n \in \mathbb{N}$. As in the proof of Lemma 2.5, we can assume that $u_n \rightarrow u$ and $v_n \rightarrow v$, for some $u,v \in \mathbb{C}$ with $|u|=|v|=1$. Furthermore, the proof of Lemma 2.5 guarantees that for some $r_1$, $r_2>0$ the open half-lines $L_{\alpha,u}=\{\alpha+tu:t>-r_1\}$ and $L_{\beta,v}=\{\beta+tv:t>-r_2\}$ are contained in $\Omega$, where $L_{\alpha,u}$ passes through point $\alpha$ and is parallel to $u$, while $L_{\beta,v}$ passes through point $\beta$ and is parallel to $v$. A short calculation shows that:
\begin{align*}
u_n-v_n=\frac{z_n-\alpha}{|z_n-\alpha|}-\frac{z_n-\beta}{|z_n-\beta|} \rightarrow 0,
\end{align*}
thus $u=v$. Therefore, $L_{\alpha,u}$ and $L_{\beta,v}$ are parallel, as desired.\\
This means that $u$ can be thought of as a direction of the set $\Omega$. Also, if we consider the set of all open half-lines contained in $\Omega$ emanating from a fixed point $\alpha$, then this set defines, in general, an open cone contained in $\Omega$ with a certain angle $\theta$. The only exception is the case of an open strip, where the above set contains only two opposite half-lines. Clearly, $\theta$ is independent of the choice of the starting point $\alpha$. If $\theta$ is greater than $180^{\circ}$, then $\Omega=\mathbb{C}$ and if $\theta=180^{\circ}$, then $\Omega$ is a open half-plane. If $\theta$ is less than $180^{\circ}$, not much can be said about the geometry of $\Omega$; it could be an open strip, the interior of an angle, the interior of a parabola or the interior of the branch of a hyperbola etc.\\\\
\textbf{Remark 2.7:} The inverse of Lemma 2.5 does not hold, an open set $\Omega \subseteq \mathbb{C}$ which is the union of open (parallel) half-lines is certainly unbounded, but not necessarily convex, even if $\Omega$ is assumed to be connected; this fact is illustrated by picking $\Omega=\{z\in \mathbb{C}:|z|>1\}$.\\
Furthermore, Lemma 2.5 and Remark 2.6 are also valid for any open, unbounded, convex set $\Omega \subseteq \mathbb{R}^d$ or $\mathbb{C}^d$, $d \in \mathbb{N}$, since the unit sphere is compact in each case. However, the compactness of the unit sphere, which is equivalent to the finiteness of the dimension of the ambient space, is essential for Lemma 2.5 and Remark 2.6 to hold, since neither of these statements are true in general, if such an $\Omega$ is a subset of an infinite-dimensional topological vector space. We will neither concern ourselves with these type of results in the present paper, nor will we provide counter-examples, though they will be investigated in future work.\\\\
\textbf{Remark 2.8:} If $\Omega \subseteq \mathbb{C}$ is a closed, unbounded, convex set, then it is the union of closed half-lines, since by repeating the proof of Lemma 2.5 and replacing $2t$ with $t$, we have that $w_n \in \Omega$ for large enough $n$ and that $w_n \rightarrow \alpha +th$, thus $\alpha +th \in \Omega$ from the closedness of $\Omega$. Clearly, the last step of that proof cannot be repeated when dealing with boundary points of $\Omega$, though it can be repeated for its interior points. Thus, for any point $\alpha \in \Omega$ we can find a closed half-line entirely contained in $\Omega$ containing point $\alpha$ and if $\alpha$ belongs to the interior $\Omega^o$ of $\Omega$, then the half-line can be chosen to be open. However, if $\Omega$ is assumed neither open, nor closed, then Lemma 2.5 does not hold in general; for instance, let $\Omega$ be the strip $\{z\in \mathbb{C}:0<$ Im$z<1 \} \cup \{0\}$ which has $0$ as a boundary point, but no half-line entirely contained in $\Omega$ exists containing $0$. In light of this remark, one can easily formulate results analogous to Lemma 2.5 and Remarks 2.6 and 2.7 for closed, unbounded, convex sets.\\\\

\section{The spaces $H_{F}^{\infty}(\Omega)$ and $\mathcal{A}_{F}(\Omega)$}
If $\Omega \subseteq \mathbb{C}$ is a domain, we denote by $H^{\infty}(\Omega)$ the space of bounded holomorphic functions on $\Omega$, endowed with the supremum norm; it is a Banach space. If $F\subseteq \mathbb{N}_0=\{0,1,...\}$ is a non-empty set, then we wish to consider the space $H_{F}^{\infty}(\Omega)$, containing all holomorphic functions $f$ on $\Omega$, such that the derivative $f^{(l)}$ belongs to $H^{\infty}(\Omega)$, for all $l\in F$. Namely:
\begin{align*}
H_{F}^{\infty}(\Omega)=\{f\in H(\Omega): f^{(l)}\in H^{\infty}(\Omega),\ \text{for all} \ l \in F\}.
\end{align*}
We define a natural topology on this space via the seminorms:
\begin{align*}
\underset{z\in \Omega}{\sup} \left | f^{(l)}(z) \right |, \ \text{for} \ l\in F \ \text{and} \ \left | f^{(l)}(z_o) \right |, \ \text{for} \ 0\leqslant l < \min F,
\end{align*}
where $z_o$ is an arbirtrary, yet fixed point in $\Omega$. We will show that $H_{F}^{\infty}(\Omega)$ is a complete metric space, hence a Fréchet space. In fact, if $F$ is finite, it is a Banach space. In any case, Baire's theorem is at our disposal.\\\\
\textbf{Theorem 3.1:} Let $\Omega \subseteq \mathbb{C}$ be a domain and $F\subseteq \mathbb{N}_0$ be a non-empty set. Then, the space $H_{F}^{\infty}(\Omega)$ with its natural topology is a complete metric space.\\\\
For the proof we need the following propositions:\\\\
\textbf{Proposition 3.2:} Let $\Omega \subseteq \mathbb{C}$ be a domain and $z_o$ be a fixed point in $\Omega$. Also, let $f_n,f$, $n\in \mathbb{N}$, be holomorphic functions on $\Omega$. Assume that $f_n' \rightarrow f'$ uniformly on compact subsets of $\Omega$ and that $f_n(z_o) \rightarrow f(z_o)$. Then, $f_n \rightarrow f$ uniformly on compact subsets of $\Omega$.\\\\
\textit{Proof:} Let $D(z,r)$ be a disc, such that $\overline{D(z,r)}\subseteq \Omega$. If $f_n(z) \rightarrow f(z)$, then by writing:
\begin{align*}
f(w)=f(z)+{\displaystyle \int_{[z,w]}f'(\zeta)d\zeta},
\end{align*}
for all $w\in\overline{D(z,r)}$, one can easily show that $f_n \rightarrow f$ uniformly on $\overline{D(z,r)}$. Therefore, the set $G=\{z\in \Omega:f_n(z) \rightarrow f(z)\}$ can easily seen to be open and closed in $\Omega$. Since $z_o \in G \neq \varnothing$ and $\Omega$ is connected, it follows that $G=\Omega$ and thus the convergence $f_n \rightarrow f$ is uniform on every closed disc contained in $\Omega$. Since every compact subset of $\Omega$ can be covered by a finite union of such discs, we conclude that $f_n \rightarrow f$ uniformly on compact subsets of $\Omega$. $\qed$\\\\
\textbf{Proposition 3.3:} Let $\Omega \subseteq \mathbb{C}$ be a domain and $z_o$ be a fixed point in $\Omega$. Also, let $(f_n)_n$ be a sequence of holomorphic functions on $\Omega$. Assume that the sequence $(f_n')_n$ is uniformly Cauchy on compact subsets of $\Omega$ and that the sequence $(f_n(z_o))_n$ is Cauchy. Then, there exists a holomorphic function $f$ on $\Omega$, such that $f_n \rightarrow f$ uniformly on compact subsets of $\Omega$. \nopagebreak \\\\
\textit{Proof:} Let $g(z)= \lim f_n'(z)$, where the convergence is uniform on every compact subset of $\Omega$. If we show that $g$ has a primitive $f$ on $\Omega$, then by adding a constant we can obtain $f(z_o)= \lim f_n(z_o)$. Then, Proposition 3.2 yields the result. Therefore, it remains to prove that $g$ has a primitive on $\Omega$, even though $\Omega$ is not assumed to be simply connected. It suffices to show that:
\begin{align*}
{\displaystyle \int_{\gamma}g(\zeta)d\zeta}=0,
\end{align*}
for all closed polygonal lines $\gamma$ in $\Omega$ ([1]). Let $\gamma$ be such a curve. Since $f_n' \rightarrow g$ uniformly on the compact set $\gamma$, it follows that:
\begin{align*}
{\displaystyle \int_{\gamma}g(\zeta)d\zeta}=\lim{\displaystyle \int_{\gamma}f_n'(\zeta)d\zeta}=0,
\end{align*}
where the last equality is true because $\gamma$ is a closed curve.$\qed$\\\\
A combination of Propositions 3.2 and 3.3 easily implies Theorem 3.1: $H_{F}^{\infty}(\Omega)$ is a Fréchet space.\\\\
We now turn our attention to the second space at hand. If $\Omega \subseteq \mathbb{C}$ is a domain, we denote by $\mathcal{A}(\Omega)$ the space of holomorphic functions on $\Omega$ possessing a continuous extension on $\overline{\Omega}$, where the closure is taken in $\mathbb{C}$. If $\Omega$ is bounded, then this space, endowed with the supremum norm, is a Banach space. If $\Omega$ is unbounded, then the topology of $\mathcal{A}(\Omega)$ is defined by the seminorms:
\begin{align*}
\underset{\underset{|z|\leqslant n}{z\in \overline{\Omega}}}{\sup} \left | f(z) \right |, \ \text{for} \ n\in \mathbb{N},
\end{align*}
and then it is a Fréchet space. If $F \subseteq \mathbb{N}_0=\{0,1,...\}$ is a non-empty set, then we wish to consider the space $\mathcal{A}_{F}(\Omega)$, containing all holomorphic functions $f$ on $\Omega$, such that the derivative $f^{(l)}$ belongs to $\mathcal{A}(\Omega)$, for all $l\in F$. Namely:
\begin{align*}
\mathcal{A}_{F}(\Omega)=\{f\in H(\Omega): f^{(l)}\in \mathcal{A}(\Omega), \ \text{for all} \ l \in F\}.
\end{align*}
The natural topology of $\mathcal{A}_{F}(\Omega)$ is the one defined by the seminorms:
\begin{align*}
\underset{\underset{|z|\leqslant n}{z\in \overline{\Omega}}}{\sup} \left | f^{(l)}(z) \right |, \ \text{for} \ l\in F, \ n\in \mathbb{N} \quad \text{and} \quad \left | f^{(l)}(z_o) \right |, \ \text{for} \ 0\leqslant l < \min F,
\end{align*}
where $z_o$ is an arbitrary, yet fixed point in $\Omega$. We will show that $\mathcal{A}_{F}(\Omega)$ is a complete metric space, hence a Fréchet space. In fact, if $\Omega$ is bounded and $F$ is finite, it is a Banach space. Thus, Baire's Theorem can be applied once again.\\\\
\textbf{Theorem 3.4:}  Let $\Omega$ be a domain and $F\subseteq \mathbb{N}_0$ be a non-empty set. Then, the space $\mathcal{A}_{F}(\Omega)$ with its natural topology is a complete metric space.\\\\
The proof of Theorem 3.4 is similar to that of Theorem 3.1, therefore, it is ommited.\\\\
We proceed to the setting of open sets. Let $\Omega$ be an open subset of $\mathbb{C}$. Then, $\Omega$ has countable connected components $\Omega_i$, $i\in I$, where $I$ is either finite, or $I=\mathbb{N}$. For every $i\in I$, we fix a point $z_i$ in $\Omega_i$. Let $F\subseteq \mathbb{N}_0=\{0,1,...\}$ be a non-empty set. We consider the spaces:
\begin{align*}
H_{F}^{\infty}(\Omega)=\{f\in H(\Omega): f^{(l)}\in H^{\infty}(\Omega),\ \text{for all} \ l \in F\}
\end{align*}
and
\begin{align*}
\mathcal{A}_{F}(\Omega)=\{f\in H(\Omega): f^{(l)}\in \mathcal{A}(\Omega), \ \text{for all} \ l \in F\}.
\end{align*}
The topology of $H_{F}^{\infty}(\Omega)$ is induced by the seminorms:
\begin{align*}
\underset{z\in \Omega}{\sup} \left | f^{(l)}(z) \right |, \ \text{for} \ l\in F \quad \text{and} \quad \left | f^{(l)}(z_i) \right |, \ \text{for} \ 0\leqslant l < \min F, \ i\in I.
\end{align*}
The topology of $\mathcal{A}_{F}(\Omega)$ is induced by the seminorms:
\begin{align*}
\underset{\underset{|z|\leqslant n}{z\in \overline{\Omega}}}{\sup} \left | f^{(l)}(z) \right |, \ \text{for} \ l\in F, \ n\in \mathbb{N} \quad \text{and} \quad \left | f^{(l)}(z_i) \right |, \ \text{for} \ 0\leqslant l < \min F, \ i\in I.
\end{align*}
By applying the previous results of this section regarding domains to each connected component of the open set $\Omega$, we deduce that $H_{F}^{\infty}(\Omega)$ and $\mathcal{A}_{F}(\Omega)$ are Fréchet spaces; the proofs of these assertions are similar to the ones in the case where $\Omega$ was a domain, only with some minor modifications. Therefore, Theorems 3.1 and 3.4 extend to the case of non-connected open sets $\Omega$. Also, Baire's theorem can be applied.\\\\
\textbf{Theorem 3.5:} Let $\Omega \subseteq \mathbb{C}$ be an open set and $F \subseteq \mathbb{N}_0$ be a non-empty set. Then, the space $H_{F}^{\infty}(\Omega)$ with its natural topology is a complete metric space.\\\\
\textbf{Theorem 3.6:} Let $\Omega \subseteq \mathbb{C}$ be an open set and $F \subseteq \mathbb{N}_0$ be a non-empty set. Then, the space $\mathcal{A}_{F}(\Omega)$ with its natural topology is a complete metric space.\\\\

\section{Two special cases of open sets and bounded functions}
In this section, we consider two special cases of open sets $\Omega \subseteq \mathbb{C}$, where $\Omega$ is:\\\\
\textbf{A:} the union of open half-lines\\
\textbf{B:} a bounded convex domain\\\\
and we prove some results regarding the space $H_{F}^{\infty}(\Omega)$, where $F \subseteq \mathbb{N}_0=\{0,1,...\}$ is a non-empty set. In particular, we examine whether the spaces $H_{F}^{\infty}(\Omega)$ and $H_{\widetilde{F}}^{\infty}(\Omega)$ coincide or not, where $\widetilde{F}=\{ l\in \mathbb{N}_0:\min F\leqslant l \leqslant \sup F\}$. In case \textbf{A}, this is true. In case \textbf{B}, we can replace $\widetilde{F}$ with the set $\widetilde{F}_0= \{  l\in \mathbb{N}_0:0\leqslant l \leqslant \sup F\}$ and then of course $H_{F}^{\infty}(\Omega)=H_{\widetilde{F}}^{\infty}(\Omega)$.\\\\
\textbf{Case A:} We present the following interpolation inequality involving derivatives of functions of one real variable, which gave us the motivation for case \textbf{A}; we refer to [3]:\\\\
\textbf{Theorem 4.1:} (Landau-Kolmogorov inequality) Let $f$ be a real- or complex-valued function defined on $I$, where $I=\mathbb{R}$ or $I=(0,+\infty)$. Assume that $f$ is $n$-times differentiable on $I$ and let:
\begin{align*}
M_k=\underset{x\in I}{\sup} \left | f^{(k)}(x) \right |,
\end{align*}
for $k=0,1,...,n$. Then, if both $M_0$ and $M_n$ are finite, the following bounds are valid:
\begin{align*}
M_k \leqslant C(n,k,I)\cdot M_0^{1-k/n}\cdot M_n^{k/n},
\end{align*}
for $k=1,...,n-1$, where $0<C(n,k,I)<+\infty$ are constants dependent only on $n,k$ and $I$.\\\\
It follows that if $f$ and $f^{(n)}$ are bounded, then all the indermediate derivatives $f^{(k)}$, $k=1,...,n-1$, are bounded as well.\\
We make some additional remarks. These constants depend on $I$ only in the following sense: $C(n,k,\mathbb{R})$ and $C(n,k,(0,+\infty))$ are different, but in each case they depend only on $n$ and $k$. Hence, from now on we will denote these constants simply by $C(n,k)$, since the dependency on $I$ is of no true significance and to emphasize the fact that we compare derivatives of order $0,k$ and $n$. Furthermore, $C(n,k)$ lie between $1$ and $\frac{\pi}{2}$ for all $n$ and $k$ and can be expressed in terms of the Favard constants, which are defined as some series of numbers; these results, among others, are due to Kolmogorov, see [3] for details and some results on the asymptotic behaviour of $C(n,k)$. Most of the literature surrounding this inequality is concerned with describing the behaviour of $C(n,k)$ and providing sharp bounds and estimates for them. We note that, for our purposes, the exact value of these constants will be irrelevant; an explicit formula is unknown anyway.\\
Also, it is clear that such an inequality holds for functions defined on any open unbounded interval, when comparing derivatives of any order. Finally, we note that we shall make use of this inequality in the following simpler form:
\begin{align*}
M_k \leqslant C(n,k)\cdot \max \{ M_0, M_n \},
\end{align*}
for $k=1,...,n-1$, which we easily derive from the previous one. Inspired by Theorem 4.1, we are now ready to prove the following:\\\\
\textbf{Theorem 4.2:} Let $\Omega \subseteq \mathbb{C}$ be an open set, which is the union of open half-lines and $F \subseteq \mathbb{N}_0$ be a non-empty set. Then, $H_{F}^{\infty}(\Omega)=H_{\widetilde{F}}^{\infty}(\Omega)$, where $\widetilde{F}=\{ l\in \mathbb{N}_0: \min F\leqslant l \leqslant \sup F\}$.\\\\
\textit{Proof:} Clearly, $H_{\widetilde{F}}^{\infty}(\Omega) \subseteq H_{F}^{\infty}(\Omega)$. For the inverse inclusion, let $f \in H_{F}^{\infty}(\Omega)$. Fix any point $p$ in $\Omega$, then there exists an open half-line $L$ contained entirely in $\Omega$, upon which lies point $p$. We can identify $L$ with the interval $I:=(-\varepsilon,+\infty)$, for any fixed $\varepsilon>0$; that is, $I$ will serve as a parametrization of $L$. Choose an $h\in \mathbb{C}$, parallel to $L$, with $|h|=1$. Consider the function $g:I \rightarrow \mathbb{C}$, defined by $g(t)=f(p+th)$. Since $f$ is holomorphic on $\Omega$, $g$ is of class $C^{\infty}$ on $I$ and $g^{(k)}(t)=f^{(k)}(p+th)\cdot h^k$, for all $t\in I$ and $k \in \mathbb{N}_0$. Therefore, $|g^{(k)}(0)|=|f^{(k)}(p)|$, for all $k \in \mathbb{N}_0$, since $|h|=1$. Observe that:
\begin{align*}
M_k:= \underset{x\in I}{\sup} \left | g^{(k)}(x) \right |=\underset{z\in L}{\sup} \left | f^{(k)}(z) \right |,
\end{align*}
for all $k \in \mathbb{N}_0$. We will show that $f \in H_{\widetilde{F}}^{\infty}(\Omega)$. To this end, we pick an $l\in \widetilde{F}\setminus F$ and we will show that $f^{(l)}\in H^{\infty}(\Omega)$ (if $l\in F$, we have nothing to prove). Next, pick any $\alpha_1, \alpha_2 \in F$ such that $\alpha_1<l<\alpha_2$. Since $f \in H_{F}^{\infty}(\Omega)$ and $\alpha_1,\alpha_2 \in F$, we have that $f^{(\alpha_1)}$ and $f^{(\alpha_2)}$ are bounded on $\Omega$, thus they are also bounded on $L\subseteq \Omega$. Consequently, $g^{(\alpha_1)}$ and $g^{(\alpha_2)}$ are bounded on $I$, that is $M_{\alpha_1}<+\infty$ and $M_{\alpha_2}<+\infty$. By invoking Theorem 4.1, we obtain the existence of some constants $C(\alpha_1,k,\alpha_2)$ satisfying:
\begin{align*}
M_k \leqslant C(\alpha_1,k,\alpha_2)\cdot \max \{ M_{\alpha_1}, M_{\alpha_2} \},
\end{align*}
for all $k\in \mathbb{N}$ satisfying $\alpha_1 < k <\alpha_2$. Since $\alpha_1<l<\alpha_2$, it follows that $M_l<+\infty$, which means that $g^{(l)}$ is also bounded on $I$. Hence:
\begin{align*}
\left | f^{(l)}(p) \right |=\left |g^{(l)}(0)  \right |\leqslant M_l\leqslant C(\alpha_1,l,\alpha_2)\cdot \max \{ M_{\alpha_1}, M_{\alpha_2} \}
\end{align*}
and for $j=\alpha_1,\alpha_2$ we have that:
\begin{align*}
M_{\alpha_j}=\underset{z\in L}{\sup} \left | f^{(\alpha_j)}(z) \right | \leqslant \underset{z\in \Omega}{\sup} \left | f^{(\alpha_j)}(z) \right |<+\infty.
\end{align*}
The last two inequalities combined imply that:
\begin{align*}
\left | f^{(l)}(p) \right | \leqslant C(\alpha_1,l,\alpha_2)\cdot \max \left \{ \underset{z\in \Omega}{\sup} \left | f^{(\alpha_1)}(z) \right |, \underset{z\in \Omega}{\sup} \left | f^{(\alpha_2)}(z) \right | \right \},
\end{align*}
from which we obtain the following:
\begin{align*}
\underset{z\in \Omega}{\sup} \left | f^{(l)}(z) \right |\leqslant C(\alpha_1,l,\alpha_2)\cdot \max \left \{ \underset{z\in \Omega}{\sup} \left | f^{(\alpha_1)}(z) \right |, \underset{z\in \Omega}{\sup} \left | f^{(\alpha_2)}(z) \right | \right \}<+\infty.
\end{align*}
Hence, $f^{(l)}$ is bounded on $\Omega$; that is $f^{(l)}\in H^{\infty}(\Omega)$ and $l$ was arbitrary. This implies that $f\in H_{\widetilde{F}}^{\infty}(\Omega)$, thus $H_{F}^{\infty}(\Omega) \subseteq H_{\widetilde{F}}^{\infty}(\Omega)$, completing the proof.$\qed$\\\\
We have proved that $H_{F}^{\infty}(\Omega)=H_{\widetilde{F}}^{\infty}(\Omega)$, where $\Omega  \subseteq \mathbb{C}$ is the union of open half-lines, $F \subseteq \mathbb{N}_0$ is a non-empty set and $\widetilde{F}=\{ l\in \mathbb{N}_0:\min F\leqslant l \leqslant \sup F\}$. So far, this is an equality between sets. We will show that the topologies of these spaces coincide as well.\\\\
We remind that if $\Omega \subseteq \mathbb{C}$ is an open set with connected components $\Omega_i$, $i\in I$, where $I$ is either finite, or $I=\mathbb{N}$ and $z_i$ is a fixed point in each $\Omega_i$, we topologize $H_{\widetilde{F}}^{\infty}(\Omega)$ via the seminorms:
\begin{align*}
\underset{z\in \Omega}{\sup} \left | f^{(l)}(z) \right |, \ \text{for} \ l\in F \quad \text{and} \quad \left | f^{(l)}(z_i) \right |, \ \text{for} \ 0\leqslant l < \min F, \ i\in I.
\end{align*}
Since $F \subseteq \widetilde{F}$ and $\min F=min\widetilde{F}$, the topology of $H_{\widetilde{F}}^{\infty}(\Omega)$ is defined by the same seminorms, in addition to the following ones:
\begin{align*}
\underset{z\in \Omega}{\sup} \left | f^{(l)}(z) \right |, \ \text{for} \ l\in \widetilde{F}\setminus F,
\end{align*}
for which we gave the following bounds during the proof of Theorem 4.2:
\begin{align*}
\underset{z\in \Omega}{\sup} \left | f^{(l)}(z) \right |\leqslant C(\alpha_1,l,\alpha_2)\cdot \max \left \{ \underset{z\in \Omega}{\sup} \left | f^{(\alpha_1)}(z) \right | ,\underset{z\in \Omega}{\sup} \left | f^{(\alpha_2)}(z) \right | \right \},
\end{align*}
where $l\in \widetilde{F}\setminus F$, $\alpha_1,\alpha_2 \in F$ satisfying $\alpha_1<l<\alpha_2$ and $C(\alpha_1,l,\alpha_2)$ are the constants mentioned in Theorem 4.1. It follows that these topologies are indeed the same. Therefore, we have proved the following statement: \\\\
\newpage \noindent \textbf{Proposition 4.3:} Let $\Omega \subseteq \mathbb{C}$ be an open set, which is the union of open half-lines and $F \subseteq \mathbb{N}_0$ be a non-empty set. Then, $H_{F}^{\infty}(\Omega)=H_{\widetilde{F}}^{\infty}(\Omega)$ as topological spaces, where $\widetilde{F}=\{ l\in \mathbb{N}_0:\min F\leqslant l \leqslant \sup F\}$. \\\\
\textbf{Remark 4.4:} A second proof of the equivalence of these topologies can be given using the Open Mapping theorem for Fréchet spaces.\\\\
Recall Lemma 2.5 from section 2. If $\Omega \subseteq \mathbb{C}$ is an unbounded convex domain, then $\Omega$ is the union of open (parallel) half-lines. Consequently, Theorem 4.2 and Proposition 4.3 are valid for such a domain $\Omega$. Therefore, we combine Lemma 2.5, Theorem 4.2 and Proposition 4.3 and we obtain the following:\\\\
\textbf{Theorem 4.5:} Let $\Omega \subseteq \mathbb{C}$ be an unbounded convex domain and $F \subseteq \mathbb{N}_0$ be a non-empty set. Then, $H_{F}^{\infty}(\Omega)=H_{\widetilde{F}}^{\infty}(\Omega)$, where $\widetilde{F}=\{ l\in \mathbb{N}_0:\min F\leqslant l \leqslant \sup F\}$.\\\\
\textbf{Proposition 4.6:} Let $\Omega \subseteq \mathbb{C}$ be an unbounded convex domain and $F \subseteq \mathbb{N}_0$ be a non-empty set. Then, $H_{F}^{\infty}(\Omega)=H_{\widetilde{F}}^{\infty}(\Omega)$ as topological spaces, where $\widetilde{F}=\{ l\in \mathbb{N}_0: \min F\leqslant l \leqslant \sup F\}$.\\\\
\textbf{Case B:} We begin with an elementary observation. Let $\Omega \subseteq \mathbb{C}$ be a bounded convex domain and $f$ be a bounded holomorphic function on $\Omega$. Using the convexity and boundedness of $\Omega$, it is easy to see that the primitive:
\begin{align*}
F(z)={\displaystyle \int_{[\alpha,z]}f(\zeta)d\zeta}
\end{align*}
of $f$, where $\alpha$ is an arbitrary, yet fixed point in $\Omega$, determining the path of integration, is Lipschitz continuous on $\Omega$. Thus, $F$ is uniformly continuous on $\Omega$, which implies that $F$ is continuously extendable on the compact set $\overline{\Omega}$. This in turn implies that $F$ is bounded on $\Omega$. Since every other primitive of $f$ differs from $F$ only by a constant, we deduce that every primitive of $f$ is bounded and uniformly continuous on $\Omega$. By making use of these statements, we prove the following:\\\\
\textbf{Theorem 4.7:} Let $\Omega \subseteq \mathbb{C}$ be an bounded convex domain and $F \subseteq \mathbb{N}_0$ be a non-empty set. Then, $H_{F}^{\infty}(\Omega)=H_{\widetilde{F}_0}^{\infty}(\Omega)$, where $\widetilde{F}_0=\{ l\in \mathbb{N}_0:0\leqslant l \leqslant \sup F\}$.\\\\
\textit{Proof:} Clearly, $H_{\widetilde{F}_0}^{\infty}(\Omega) \subseteq H_{F}^{\infty}(\Omega)$. For the inverse inclusion, let $f \in H_{F}^{\infty}(\Omega)$. We will show that $f \in H_{\widetilde{F}_0}^{\infty}(\Omega)$. To this end, we pick an $l\in \widetilde{F}_0\setminus F$ and we will show that $f^{(l)}\in H^{\infty}(\Omega)$ (if $l\in F$, we have nothing to prove). Next, pick any $\alpha \in F$ such that $\alpha>l$. Since $f \in H_{F}^{\infty}(\Omega)$ and $\alpha \in F$, we have that $f^{(\alpha)}$ is bounded on $\Omega$. By integrating the bounded function $f^{(\alpha)}$ repeatedly, we deduce that all the functions $f^{(k)}$, $k=0,...,\alpha-1$, are also bounded on $\Omega$; this follows from the earlier discussion. Hence, $f^{(l)}$ is bounded on $\Omega$; that is $f^{(l)}\in H^{\infty}(\Omega)$ and $l$ was arbitrary. This implies that $f\in H_{\widetilde{F}_0}^{\infty}(\Omega)$, thus $H_{F}^{\infty}(\Omega) \subseteq H_{\widetilde{F}_0}^{\infty}(\Omega)$, completing the proof.$\qed$\\\\
We have proved that $H_{F}^{\infty}(\Omega)=H_{\widetilde{F}_0}^{\infty}(\Omega)$, where $\Omega \subseteq \mathbb{C}$ is a bounded convex domain, $F \subseteq \mathbb{N}_0$ is a non-empty set and $\widetilde{F}_0=\{ l\in \mathbb{N}_0:0\leqslant l \leqslant \sup F\}$. So far, this is an equality between sets. We will show that the topologies of these spaces coincide as well.\\\\
We remind that if $\Omega \subseteq \mathbb{C}$ is a domain and $z_o$ is a fixed point in $\Omega$, we topologize $H_{F}^{\infty}(\Omega)$ via the seminorms:
\begin{align*}
\underset{z\in \Omega}{\sup} \left | f^{(l)}(z) \right |, \ \text{for} \ l\in F \quad \text{and} \quad \left | f^{(l)}(z_o) \right |, \ \text{for} \ 0\leqslant l < \min F,
\end{align*}
while the topology of $H_{\widetilde{F}_0}^{\infty}(\Omega)$ is induced by the seminorms:
\begin{align*}
\underset{z\in \Omega}{\sup} \left | f^{(l)}(z) \right |, \ \text{for} \ l\in \widetilde{F}_0.
\end{align*}
Since $F \subseteq \widetilde{F}_0$, for any $l\in F$ the seminorm $\underset{z\in \Omega}{\sup} \left | f^{(l)}(z) \right |$ is taken into account in both of these topologies. If $0\leqslant l < \min F$, then $l\in \widetilde{F}_0$ and obviously:
\begin{align*}
\left | f^{(l)}(z_o) \right | \leqslant \underset{z\in \Omega}{\sup} \left | f^{(l)}(z) \right |.
\end{align*}
Hence, the topology of $H_{\widetilde{F}_0}^{\infty}(\Omega)$ is finer that the topology of $H_{F}^{\infty}(\Omega)$. Now pick any $l\in \widetilde{F}_0$. Then, either $l\in F$ and thus the seminorm $\underset{z\in \Omega}{\sup} \left | f^{(l)}(z) \right |$ is taken into account in both of these topologies, or $l\notin F$. In the latter case, pick any $\alpha \in F$ such that $\alpha>l$. For all $z\in \Omega$ we have that:
\begin{align*}
\left |f^{(\alpha-1)}(z)  \right | &=\left |{\displaystyle \int_{[z_o,z]}f^{(\alpha)}(\zeta)d\zeta}+f^{(\alpha-1)}(z_o)  \right | \\\\
&\leqslant \underset{w\in \Omega}{\sup} \left | f^{(\alpha)}(w) \right | \cdot |z-z_o|+\left | f^{(\alpha-1)}(z_o) \right | \\\\
&\leqslant \underset{w\in \Omega}{\sup} \left | f^{(\alpha)}(w) \right | \cdot \text{diam}(\Omega)+ \left | f^{(\alpha-1)}(z_o) \right |,
\end{align*}\\
since $|z-z_o|\leqslant$ diam$(\Omega)<+\infty$, for all $z\in \Omega$. Continuing in this manner, we have:\\
\begin{align*}
\underset{z\in \Omega}{\sup}\left | f^{(\alpha-m)}(z) \right | &\leqslant \underset{z\in \Omega}{\sup}\left | f^{(\alpha)}(z) \right|\cdot (\text{diam}(\Omega))^{m}+\left |f^{(\alpha-1)}(z_o)  \right |\cdot (\text{diam}(\Omega))^{m-1}+\\\\
&...+\left |f^{(\alpha-m+1)}(z_o)  \right |\cdot \text{diam}(\Omega)+ \left |f^{(\alpha-m)}(z_o)  \right |,
\end{align*}\\
for $m=1,...,\alpha$. Choosing $m=\alpha-l$ and writing the above inequality in a brief form, we obtain:\\
\begin{align*}
\underset{z\in \Omega}{\sup}\left | f^{(l)}(z) \right | \leqslant \underset{z\in \Omega}{\sup}\left | f^{(\alpha)}(z) \right|\cdot (\text{diam}(\Omega))^{a-l}+\sum_{k=l}^{a-1} \left | f^{(k)}(z_o) \right | \cdot (\text{diam}(\Omega))^{k-l}.
\end{align*}\\
\newpage \noindent If $(f_n)_n$ is a sequence in $H_{F}^{\infty}(\Omega)$ and $f\in H_{F}^{\infty}(\Omega)$, such that $f_n \rightarrow f$ in the topology of $H_{F}^{\infty}(\Omega)$, then by Weierstrass's theorem and Proposition 3.2, in combination with the pre- vious inequality, one can easily deduce that $f_n \rightarrow f$ in the topology of $H_{\widetilde{F}_0}^{\infty}(\Omega)$. Hence, these topologies are indeed the same. Therefore, we have proved the following statement:\\\\
\textbf{Proposition 4.8:} Let $\Omega \subseteq \mathbb{C}$ be a bounded convex domain and $F \subseteq \mathbb{N}_0$ be a non-empty set. Then, $H_{F}^{\infty}(\Omega)=H_{\widetilde{F}_0}^{\infty}(\Omega)$ as topological spaces, where $\widetilde{F}_0=\{ l\in \mathbb{N}_0:0\leqslant l \leqslant \sup F\}$.\\\\
\textbf{Remark 4.9:} A second proof of the equivalence of these topologies can be given using the Open Mapping theorem for Fréchet spaces.\\\\
\textbf{Remark 4.10:} In [11], see also [9], a Jordan domain $\Omega \subseteq \mathbb{C}$ has been constructed, suppor- ting a bounded holomorphic function $g$, so that its primitive $G$ is unbounded. Thus, for this domain $\Omega$, the spaces $H_{F}^{\infty}(\Omega)$ and $H_{\widetilde{F}_0}^{\infty}(\Omega)$ are different for some non-empty set $F \subseteq \mathbb{N}_0$, such that $0\notin F$ and $1\in F$. This is certainly true for $F=\{1\}$ and $\widetilde{F}_0=\{0,1\}$.\\\\
\textbf{Remark 4.11:} Suppose that $\Omega \subseteq \mathbb{C}$ is a simply connected domain, for which a constant $0<M<+\infty$ exists, with the property that any two points $p,q \in \Omega$ can be joined by a rectifiable curve $\gamma_{p,q}$ in $\Omega$, with length bounded by $M$; then clearly, $\Omega$ is bounded with diam$(\Omega) \leqslant M$. Then, all the results obtained in case \textbf{B} of this section for bounded convex domains are still valid for such a domain $\Omega$. Because then, if $f$ is a bounded holomorphic function on $\Omega$, its primitive is bounded by $M\cdot \underset{z\in \Omega}{\sup}\left | f(z) \right |+|f(z_o)|$, where $z_o$ is a fixed point in $\Omega$. This condition has been used in [11]. More recently, it has been proven in [12] that his condition is necessary and sufficient for a simply connected domain $\Omega$, in order for the primitive of any bounded holomorphic function on $\Omega$ to be also bounded; this condition is connected to the boundedness of the integration operator.\\\\

\section{The special case of convex domains and continuously extendable functions}
In this section, we examine if analogues of Theorems 4.5 and 4.7 are valid for the space $\mathcal{A}_{F}(\Omega)$, where $F \subseteq \mathbb{N}_0=\{0,1,...\}$ is a non-empty set and $\Omega$ is a convex domain. That is, if $\mathcal{A}_{F}(\Omega)=\mathcal{A}_{\widetilde{F}}(\Omega)$ when $\Omega$ is unbounded and $\widetilde{F}=\{ l\in \mathbb{N}_0:\min F\leqslant l \leqslant \sup F\}$ and if $\mathcal{A}_{F}(\Omega)=\mathcal{A}_{\widetilde{F}_0}(\Omega)$ when $\Omega$ is bounded and $\widetilde{F}_0=\{ l\in \mathbb{N}_0:0\leqslant l \leqslant \sup F\}$. We will show that $\mathcal{A}_{F}(\Omega)=\mathcal{A}_{\widetilde{F}_0}(\Omega)$ for any convex domain $\Omega\subseteq \mathbb{C}$, regardless of whether $\Omega$ is bounded or unbounded. Of course, this implies that $\mathcal{A}_{F}(\Omega)=\mathcal{A}_{\widetilde{F}}(\Omega)$.\\\\
\textbf{Theorem 5.1:} Let $\Omega \subseteq \mathbb{C}$ be a convex domain and $F \subseteq \mathbb{N}_0$ be a non-empty set. Then, $\mathcal{A}_{F}(\Omega)=\mathcal{A}_{\widetilde{F}_0}(\Omega)$, where $\widetilde{F}_0=\{ l\in \mathbb{N}_0:0\leqslant l \leqslant \sup F\}$.\\\\
\textit{Proof:} Clearly $\mathcal{A}_{\widetilde{F}_0}(\Omega) \subseteq \mathcal{A}_{F}(\Omega)$. For the inverse inclusion, let $f\in \mathcal{A}_{F}(\Omega)$. We will show that $f\in \mathcal{A}_{\widetilde{F}_0}(\Omega)$. To this end, we pick an $l\in \widetilde{F}_0\setminus F$ and we will show that $f^{(l)}\in \mathcal{A}(\Omega)$ (if $l\in F$, we have nothing to prove). Next, pick any $\alpha \in F$ such that $\alpha > l$. Since $f\in \mathcal{A}_{F}(\Omega)$ and $\alpha \in F$, we have that $f^{(\alpha)}$ is continuously extendable on $\overline{\Omega}$. If $\Omega$ is bounded, then $f^{(\alpha)}$ is bounded on $\Omega$ and by integrating $f^{(\alpha)}$ repeatedly we deduce that $f^{(l)}$ is uniformly continuous on $\Omega$. It follows that $f^{(l)}$ is continuously extendable on $\overline{\Omega}$. If $\Omega$ is unbounded, we work in a similar way, but only locally. Fix any point $\zeta \in \partial \Omega$ and consider the sets $\Omega_m=\Omega \cap D(0,m)$, for $m\in \mathbb{N}$. Then, $\zeta$ lies within some set $\overline{\Omega_m}^o$, where the interior is relative to $\overline{\Omega}$. Indeed, for some $k\in \mathbb{N}$ satisfying $|\zeta| \leqslant k$, one can easily see that:
\begin{align*}
\zeta \in \overline{\Omega} \cap \overline{D(0,k)} \subseteq \left ( \overline{\Omega \cap D(0,k+1)} \right )^o \subseteq \overline{\Omega \cap D(0,k+1)}.
\end{align*}
Since $f^{(\alpha)}$ is continuously extendable on $\overline{\Omega}$, it is also continuously extendable on the compact set $\overline{\Omega_{k+1}}$. Hence, $f^{(\alpha)}$ is bounded on the bounded convex domain $\Omega_{k+1}$ and by integrating $f^{(\alpha)}$ repeatedly we deduce that $f^{(l)}$ is uniformly continuous on $\Omega_{k+1}$. It follows that $f^{(l)}$ is continuously extendable on $\overline {\Omega_{k+1}}$ and thus extends continuously at point $\zeta$ which was arbitrary. It follows that $f^{(l)}$ is continuously extendable on $\overline{\Omega}$. Hence, $f^{(l)}$ is continuously extendable on $\overline{\Omega}$, whether $\Omega$ is bounded or unbounded; that is $f^{(l)}\in \mathcal{A}(\Omega)$ and $l$ was arbitrary. This implies that $f\in \mathcal{A}_{\widetilde{F}_0}(\Omega)$, thus $\mathcal{A}_{F}(\Omega) \subseteq \mathcal{A}_{\widetilde{F}_0}(\Omega)$, completing the proof.$\qed$\\\\
We have proved that $\mathcal{A}_{F}(\Omega)=\mathcal{A}_{\widetilde{F}_0}(\Omega)$, where $\Omega \subseteq \mathbb{C}$ is a convex domain, $F \subseteq \mathbb{N}_0$ is a non-empty set and $\widetilde{F}_0=\{ l\in \mathbb{N}_0:0\leqslant l \leqslant \sup F\}$. So far, this is an equality between sets. We will show that the topologies of these spaces coincide as well.\\\\
We remind that if $\Omega \subseteq \mathbb{C}$ is a domain and $z_o$ is a fixed point in $\Omega$, we topologize $\mathcal{A}_{F}(\Omega)$ via the seminorms:
\begin{align*}
\underset{\underset{|z|\leqslant n}{z\in \overline{\Omega}}}{\sup} \left | f^{(l)}(z) \right |, \ \text{for} \ l\in F, \ n\in \mathbb{N} \quad \text{and} \quad \left | f^{(l)}(z_o) \right |, \ \text{for} \ 0\leqslant l < \min F,
\end{align*}
while the topology of $\mathcal{A}_{\widetilde{F}_0}(\Omega)$ is induced by the seminorms:
\begin{align*}
\underset{\underset{|z|\leqslant n}{z\in \overline{\Omega}}}{\sup} \left | f^{(l)}(z) \right |, \ \text{for} \ l\in \widetilde{F}_0, \ n\in \mathbb{N}.
\end{align*}
Since $F \subseteq \widetilde{F}_0$, for any $l\in F$ the seminorms $\underset{\underset{|z|\leqslant n}{z\in \overline{\Omega}}}{\sup} \left | f^{(l)}(z) \right |$, $n\in \mathbb{N}$, are taken into account in both of these topologies. If $0\leqslant l < \min F$, then $l\in \widetilde{F}_0$ and obviously:
\begin{align*}
\left | f^{(l)}(z_o) \right | \leqslant \underset{\underset{|z|\leqslant n_o}{z\in \overline{\Omega}}}{\sup} \left | f^{(l)}(z) \right |,
\end{align*}
for some $n_o\in \mathbb{N}$ satisfying $|z_o| \leqslant n_o$. Hence, the topology of $\mathcal{A}_{\widetilde{F}_0}(\Omega)$ is finer that the topology of $\mathcal{A}_{F}(\Omega)$. Now pick any $l\in \widetilde{F}_0$. Then, either $l\in F$ and thus the seminorms $\underset{\underset{|z| \leqslant n}{z\in \overline{\Omega}}}{\sup} \left | f^{(l)}(z) \right |$, $n\in \mathbb{N}$, are taken into account in both of these topologies, or $l\notin F$. In the latter case, pick any $\alpha \in F$ such that $\alpha>l$. If $\Omega_m=\Omega \cap D(0,m)$, for $m\in \mathbb{N}$, then notice that:
\begin{align*}
\{z\in \overline{\Omega}: |z|\leqslant m\}=\overline{\Omega} \cap \overline{D(0,m)} &\subseteq \overline{\Omega \cap D(0,m+1)} = \overline{\Omega_{m+1}}\\
&\subseteq \overline{\Omega} \cap \overline{D(0,m+1)}=\{z\in \overline{\Omega}: |z|\leqslant m+1\},
\end{align*}
for all $m \in \mathbb{N}$. From this fact and by applying the inequalities obtained during the proof of Proposition 4.8 for the closure of the bounded convex domains $\Omega_{m}$, $m\in \mathbb{N}$, (the continuity of functions in $\mathcal{A}(\Omega)$ on every $\overline{\Omega_{m}}$ guarantees that taking supremum over $\Omega_{m}$ or $\overline{\Omega_{m}}$ is the same and finite in each case) we have that:
\begin{align*}
&\qquad\qquad \underset{\underset{|z|\leqslant m}{z\in \overline{\Omega}}}{\sup} \left | f^{(l)}(z) \right | \leqslant \underset{z\in \overline{\Omega_{m+1}}}{\sup}\left | f^{(l)}(z) \right | \\
&\qquad\leqslant \underset{z\in \overline{\Omega_{m+1}}}{\sup}\left | f^{(\alpha)}(z) \right|\cdot (\text{diam}(\Omega_{m+1}))^{a-l}+\sum_{k=l}^{a-1} \left | f^{(k)}(z_o) \right | \cdot (\text{diam}(\Omega_{m+1}))^{k-l}\\
&\qquad\leqslant\underset{\underset{|z|\leqslant m+1}{z\in \overline{\Omega}}}{\sup} \left | f^{(\alpha)}(z) \right | \cdot (\text{diam}(\Omega_{m+1}))^{\alpha-l}+\sum_{k=l}^{a-1} \left | f^{(k)}(z_o) \right | \cdot (\text{diam}(\Omega_{m+1}))^{k-l},
\end{align*}\\
where diam$(\Omega_m)<+\infty$, for all $m\in \mathbb{N}$. If $(f_n)_n$ is a sequence in $\mathcal{A}_{F}(\Omega)$ and $f\in \mathcal{A}_{F}(\Omega)$, such that $f_n \rightarrow f$ in the topology of $\mathcal{A}_{F}(\Omega)$, then by Weierstrass's theorem and Proposition 3.2, in combination with the previous inequality, one can easily deduce that $f_n \rightarrow f$ in the topology of $\mathcal{A}_{F}(\Omega)$. Hence, these topologies are indeed the same. Therefore, we have proved the following statement:\\\\
\textbf{Proposition 5.2:} Let $\Omega \subseteq \mathbb{C}$ be a convex domain and $F \subseteq \mathbb{N}_0$ be a non-empty set. Then, $\mathcal{A}_{F}(\Omega)=\mathcal{A}_{\widetilde{F}_0}(\Omega)$ as topological spaces, where $\widetilde{F}_0=\{ l\in \mathbb{N}_0:0\leqslant l \leqslant \sup F\}$.\\\\
\textbf{Remark 5.3:} A second proof of the equivalence of these topologies can be given using the Open Mapping theorem for Fréchet spaces.\\\\
\textbf{Remark 5.4} For the Jordan domain $\Omega \subseteq \mathbb{C}$ mentioned in Remark 4.10, we have that the function $g$ constructed in [11] is continuously extendable on $\overline{\Omega}$, but its primitive $G$ is not, since $G$ is unbounded on $\Omega$ and $\overline{\Omega}$ is compact. Thus, for this domain $\Omega$, the spaces $\mathcal{A}_{F}(\Omega)$ and $\mathcal{A}_{\widetilde{F}_0}(\Omega)$ are different for some non-empty set $F\subseteq \mathbb{N}_0$, such that $0\notin F$ and $1\in F$. This is certainly true for $F=\{1\}$ and $\widetilde{F}_0=\{0,1\}$.\\\\

\section{Non-extendability in $H_{F}^{\infty}(\Omega)$ and $\mathcal{A}_{F}(\Omega)$}
In this section, we deal with the notion of non-extendability of holomorphic functions in the spaces $H_{F}^{\infty}(\Omega)$ and $\mathcal{A}_{F}(\Omega)$, where $\Omega \subseteq \mathbb{C}$ is an open set and $F\subseteq \mathbb{N}_0=\{0,1,...\}$ is a non-empty set. An immediate corollary of Proposition 3.2 is the following:\\\\
\textbf{Proposition 6.1:} Let $\Omega \subseteq \mathbb{C}$ be an open set and $F\subseteq \mathbb{N}_0$ be a non-empty set. Also, let $f_n,f$, $n\in \mathbb{N}$, be holomorphic functions on $\Omega$. If either:\\
(i) $f_n,f \in H_{F}^{\infty}(\Omega)$ for $n\in \mathbb{N}$ and $f_n \rightarrow f$ in the topology of $H_{F}^{\infty}(\Omega)$, or\\
(ii) $f_n,f \in \mathcal{A}_{F}(\Omega)$ for $n\in \mathbb{N}$ and $f_n \rightarrow f$ in the topology of $\mathcal{A}_{F}(\Omega)$,\\
then $f_n \rightarrow f$ uniformly on compact subsets of $\Omega$, therefore $f_n \rightarrow f$ pointwise.\\\\
This enables us to prove the following generic results:\\\\
\textbf{Theorem 6.2:} Let $\Omega \subseteq \mathbb{C}$ be an open set and $F\subseteq \mathbb{N}_0$ be a non-empty set. Then, the set $S_{H_{F}^{\infty}(\Omega)}$ of functions in $H_{F}^{\infty}(\Omega)$ which are non-extendable is either void, or a dense and $G_\delta$ subset of $H_{F}^{\infty}(\Omega)$.\\\\
\textit{Proof:} Assume that $S_{H_{F}^{\infty}(\Omega)} \neq \varnothing$. Then, by combining the completeness of the metric space $H_{F}^{\infty}(\Omega)$ with condition (i) of Proposition 6.1, we deduce that the assumptions of Theorem 2.3 for $X(\Omega)=H_{F}^{\infty}(\Omega)\subseteq H(\Omega)$ are verified. Thus, $S_{H_{F}^{\infty}(\Omega)}$ is dense and $G_\delta$ in $H_{F}^{\infty}(\Omega)$.$\qed$\\\\
\textbf{Theorem 6.3:} Let $\Omega \subseteq \mathbb{C}$ be an open set and $F\subseteq \mathbb{N}_0$ be a non-empty set. Then, the set $S_{\mathcal{A}_{F}(\Omega)}$ of functions in $\mathcal{A}_{F}(\Omega)$ which are non-extendable is either void, or a dense and $G_\delta$ subset of $\mathcal{A}_{F}(\Omega)$.\\\\
\textit{Proof:} Similar to the proof of Theorem 6.2. Assume that $S_{\mathcal{A}_{F}(\Omega)} \neq \varnothing$. Then, by combining the completeness of the metric space $\mathcal{A}_{F}(\Omega)$ with condition (ii) of Proposition 6.1, we deduce that the assumptions of Theorem 2.3 for $X(\Omega)=\mathcal{A}_{F}(\Omega)\subseteq H(\Omega)$ are verified. Thus, $S_{\mathcal{A}_{F}(\Omega)}$ is dense and $G_\delta$ in $\mathcal{A}_{F}(\Omega)$.$\qed$\\\\
Next, we give examples and investigate whether $S_{H_{F}^{\infty}(\Omega)}$ and $S_{\mathcal{A}_{F}(\Omega)}$ are empty or not:\\\\
\textbf{Example 6.4:} Let $U\subseteq \mathbb{C}$ be a domain and $K\subseteq U$ be a compact set which is remo- vable for bounded holomorphic functions; that is, its analytic capacity $\gamma(K)$ is zero. For instance, $K$ could be a singleton or a planar Cantor-type set, obtained by removing corner quarters; see [8]. Let $\Omega=U\setminus K$. Then, it is easy to see that every $f\in H_{F}^{\infty}(\Omega)$ is extendable for any choice of $F$, provided that $0\in F$. Thus, $S_{H_{F}^{\infty}(\Omega)}=\varnothing$.\\\\
\textbf{Example 6.5:} Let $U\subseteq \mathbb{C}$ be a domain and $K\subseteq U$ be a compact set with continuous analytic capacity $\alpha(K)$ equal to zero. Let $\Omega=U\setminus K$. Then, it is easy to see that every $f\in \mathcal{A}_{F}(\Omega)$ is extendable for any choice of $F$, provided that $0\in F$. Thus, $S_{\mathcal{A}_{F}(\Omega)}=\varnothing$.\\
If we use a result from [7], [13], we see that the same holds if $K$ is any closed subset of $\mathbb{C}$ with $\alpha(K)=0$; for instance, $K$ could be a straight line, a line segment, a circular arc, a circle, an analytic curve or the boundary of a convex set.\\\\
\textbf{Example 6.6:} Let $U\subseteq \mathbb{C}$ be a domain and $K\subseteq U$ be a singleton, or more generally, a compact set containing an isolated point. Let $\Omega=U\setminus K$. Then, every holomorphic function $f$ which belongs to $H_{F}^{\infty}(\Omega)$ or $\mathcal{A}_{F}(\Omega)$ is extendable, for any choice of $F$. Thus, $S_{H_{F}^{\infty}(\Omega)}=S_{\mathcal{A}_{F}(\Omega)}=\varnothing$.\\
Indeed, let $\alpha=\min F$ and $\zeta$ be an isolated point of $K$. If $f \in H_{F}^{\infty}(\Omega)$, then by Riemann's theorem on removable singularities, $f^{(\alpha)}$ is holomorphic on a sufficiently small disc $D(\zeta,r)$ contained in $\Omega$, for some $r>0$. If $f \in \mathcal{A}_{F}(\Omega)$, then $f^{(\alpha)}$ is holomorphic on $D(\zeta,r)$ by definition. In any case, since this disc is a bounded convex domain, by integrating $f^{(\alpha)}$ repeatedly we conclude that $f$ is extendable on $D(\zeta,r)$. Hence $S_{H_{F}^{\infty}(\Omega)}=S_{\mathcal{A}_{F}(\Omega)}=\varnothing$.\\\\
\textbf{Example 6.7:} Let $\Omega \subseteq \mathbb{C}$ be a domain, such that every point $\zeta \in \partial \Omega$ is the limit of a sequence $(z_n)_n$ of points contained in $\left ( \overline {\Omega} \right )^c$. Then, for any choice of $F$, the sets $S_{H_{F}^{\infty}(\Omega)}$ and $S_{\mathcal{A}_{F}(\Omega)}$ are dense and $G_\delta$ in $H_{F}^{\infty}(\Omega)$ and $\mathcal{A}_{F}(\Omega)$, respectively.\\
We will use a result from [14] regarding non-extendability, which was stated in Theorem 2.3. Pick any two discs $D_1$ and $D_2$ as in Definition 2.1 and a point $\zeta \in \partial \Omega \cap D_2$. By our assumption, there exists a point $w\in D_2 \setminus \overline {\Omega}$. Consider the function $f(z)=\frac{1}{z-w}$, which belongs to $H_{F}^{\infty}(\Omega)\cap \mathcal{A}_{F}(\Omega)$. Since this function restricted to $D_1$ is equal to $f|_{D_1}(z)=\frac{1}{z-w}$, by analytic continuation we have that its only holomorphic extension on $D_2 \setminus \{w\}$ is the function $g(z)=\frac{1}{z-w}$, which has a pole at $w\in D_2$. Thus, $f$ does not possess a holomorphic extension on $D_2$. This means that $f$ is non-extendable; that is $f\in S_{H_{F}^{\infty}(\Omega)}\cap S_{\mathcal{A}_{F}(\Omega)}$. Hence, $S_{H_{F}^{\infty}(\Omega)}\neq \varnothing$ and $S_{\mathcal{A}_{F}(\Omega)}\neq \varnothing$, therefore the sets $S_{H_{F}^{\infty}(\Omega)}$ and $S_{\mathcal{A}_{F}(\Omega)}$ are dense and $G_\delta$ in $H_{F}^{\infty}(\Omega)$ and $\mathcal{A}_{F}(\Omega)$, respectively, as Theorems 6.2 and 6.3 indicate.\\\\
\textbf{Example 6.8:} Let $\Omega \subseteq \mathbb{C}$ be a domain bounded by a finite set of disjoint Jordan curves. Then, for any choice of $F$, the sets $S_{H_{F}^{\infty}(\Omega)}$ and $S_{\mathcal{A}_{F}(\Omega)}$ are dense and $G_\delta$ in $H_{F}^{\infty}(\Omega)$ and $\mathcal{A}_{F}(\Omega)$, respectively. Clearly, this is a particular case of Example 6.7.\\\\

\section{Two more dichotomy results}
In this section, we prove two more dichotomy results regarding boundedness or unboundedness of derivatives of functions in the spaces $H_{F}^{\infty}(\Omega)$ and $\mathcal{A}_{F}(\Omega)$, where $\Omega \subseteq \mathbb{C}$ is an open set and $F \subseteq \mathbb{N}_0=\{0,1,...\}$ is a non-empty set. We will use the following result from [15]:\\\\
\textbf{Proposition 7.1:} Let $V$ be a topological vector space over the field $\mathbb{R}$ or $\mathbb{C}$ and let $X$ be a non-empty set. Denote by $\mathbb{C}^X$ the set of all complex valued functions on $X$ and consider a linear operator $T:V \rightarrow \mathbb{C}^X$ with the property that the mapping $V\ni \alpha \mapsto T_{x}(a)=T(\alpha)(x)\in \mathbb{C}$ is continuous, for all $x\in X$; observe that this assumption is weaker than $T$ being continuous. Let $S=S(T,V,X)=\{\alpha \in V: T(\alpha)$ is unbounded on $X\}$. Then, either $S=\varnothing$, or $S$ is a dense and $G_\delta$ subset of $V$.\\\\
Note that in Proposition 7.1., the space $V$ is not assumed to be a complete metric space. \\\\
Let $\Omega \subseteq \mathbb{C}$ be an open set and let $V$ be one of the topological vector spaces $H_{F}^{\infty}(\Omega)$ or $\mathcal{A}_{F}(\Omega)$, endowed with its natural topology, where $F$ is a non-empty subset of $\mathbb{N}_0$. Let $X$ be any subset of $\Omega$ and $l\in \mathbb{N}_0$. Then, the function $V\ni f \mapsto T_l(f)(z)=f^{(l)}(z) \in \mathbb{C}$ is conti- nuous, for all $z\in X$; this follows from Weierstrass's theorem if $l \geqslant \min F$ and from Proposition 3.2 if $0 \leqslant l < \min F$. Thus, the corresponding set $S=S_l$ is either empty, or dense and $G_\delta$ in the space $V$. In particular, the above holds true for $V=H_{F}^{\infty}(\Omega)$ and $X=\Omega$. Thus, we have proved the following:\\\\
\textbf{Theorem 7.2:} Let $\Omega \subseteq \mathbb{C}$ be an open set, $F \subseteq \mathbb{N}_0$ be a non-empty set and $l \in \mathbb{N}_0$. Then, either for every $f \in H_{F}^{\infty}(\Omega)$ the derivative $f^{(l)}$ is bounded on $\Omega$, or generically for every $f \in H_{F}^{\infty}(\Omega)$ the derivative $f^{(l)}$ is unbounded on $\Omega$.\\\\
If $l\in F$, then obviously for every $f \in H_{F}^{\infty}(\Omega)$ the derivative $f^{(l)}$ is bounded on $\Omega$. If $\Omega$ is a bounded convex domain and $l \leqslant \sup F$, then for every $f \in H_{F}^{\infty}(\Omega)$ the derivative $f^{(l)}$ is bounded on $\Omega$, according to Theorem 4.7.\\\\
In [11], see also [9], a Jordan domain $\Omega$ was constructed, such that a function $g:\overline{\Omega} \rightarrow \mathbb{C}$ continuous on $\overline{\Omega}$ and holomorphic on $\Omega$ has an unbounded primitive on $\Omega$. Let us call this primitive $G$; then $G \in H_{\{1\}}^{\infty}(\Omega)$, but for $l=0$ the function $G^{(0)}=G$ is unbounded on $\Omega$. Thus, in this domain $\Omega$, generically every function $f \in H_{\{1\}}^{\infty}(\Omega)$ has the property that $f^{(0)}=f$ is unbounded on $\Omega$. It follows that $H_{\{0,1\}}^{\infty}(\Omega)$ is meager in $H_{\{1\}}^{\infty}(\Omega)$ for this parti- cular domain $\Omega$. In general, either $H_{\{0,1\}}^{\infty}(\Omega)=H_{\{1\}}^{\infty}(\Omega)$, or $H_{\{0,1\}}^{\infty}(\Omega)$ is meager in $H_{\{1\}}^{\infty}(\Omega)$, for any open set $\Omega \subseteq \mathbb{C}$.\\\\
Let $\Omega=\mathbb{D}$ be the open unit disc and let $w(z)=(z-1)\cdot exp\frac{z+1}{z-1}$. Then, $w \in \mathcal{A}(\mathbb{D}) \subseteq H^{\infty}(\mathbb{D})=H_{\{0\}}^{\infty}(\mathbb{D})$ and $w'$ is unbounded on $\mathbb{D}$, thus generically every function $f \in H_{\{0\}}^{\infty}(\mathbb{D})$ has the property that $f^{(1)}=f'$ is unbounded on $\mathbb{D}$. It follows that $H_{\{0,1\}}^{\infty}(\mathbb{D})$ is meager in $H_{\{0\}}^{\infty}(\mathbb{D})$. More generally, if $F$ is finite and $l>\max F$, then generically for every $f\in H_{F}^{\infty}(\mathbb{D})$ the derivative $f^{(l)}$ is unbounded on $\mathbb{D}$ and $H_{F\cup \{l\}}^{\infty}(\mathbb{D})$ is meager in $H_{F}^{\infty}(\mathbb{D})$.\\\\
It remains open to give an example of a domain $\Omega \subseteq \mathbb{C}$, supporting a holomorphic function $f$, so that $f^{(0)}=f$  and $f^{(2)}$ are bounded on $\Omega$, but $f^{(1)}=f'$ is unbounded. We believe that such a domain $\Omega$ exists. Moreover, we think that a complete metric topology can be defined on the set of all domains (contained in the open unit disc), so that for the generic domain $\Omega$, there exists a holomorphic function $f$ on $\Omega$, such that $f$  and $f^{(2)}$ are bounded, but $f'$ is not. More generally, we think that for every non-empty set $F \subseteq \mathbb{N}_0$ and $l \notin F$, $\min F<l<\sup F$, for the generic domain $\Omega \subseteq \mathbb{C}$, there exists an $f\in H_{F}^{\infty}(\Omega)$ such that $f^{(l)}$ is unbounded and $H_{F\cup \{l\}}^{\infty}(\Omega)$ is meager in $H_{F}^{\infty}(\Omega)$. But we do not have a proof and these assertions remain open. \\\\
Next, consider the space $\mathcal{A}_{F}(\Omega)$, where $\Omega \subseteq \mathbb{C}$ is an open set and $F$ is a non-empty subset of $\mathbb{N}_0$. If $\Omega$ is bounded, then for every $f \in \mathcal{A}_{F}(\Omega)$, all derivatives $f^{(l)}$, $l \in F$, are also bounded on $\Omega$. Assume that $\Omega$ is unbounded, we will apply Proposition 7.1. We set $X=\Omega$ and let $T_l:\mathcal{A}_{F}(\Omega) \rightarrow \mathbb{C}^X$ be the function $T_l(f)(z)=f^{(l)}(z)$, $z \in X$, where $l$ is a fixed element of $\mathbb{N}_0$. Then, the assumptions of Proposition 7.1 are easily verified. Therefore, the set $S_l=\{f\in\mathcal{A}_{F}(\Omega):f^{(l)}$ is unbounded on $\Omega \}$ is either void, or dense and $G_\delta$ in $\mathcal{A}_{F}(\Omega)$. But the function $f(z)=z^{l+1}$ belongs to $\mathcal{A}_{F}(\Omega)$. Therefore, the set $S_l$ is dense and $G_\delta$ in $\mathcal{A}_{F}(\Omega)$. Baire's theorem implies that the set $S=\bigcap_{l\in \mathbb{N}_0}S_l$ is also dense and $G_\delta$ in $\mathcal{A}_{F}(\Omega)$. Thus, we have proved the following:\\\\
\textbf{Proposition 7.3:} Let $\Omega \subseteq \mathbb{C}$ be an unbounded open set and $F \subseteq \mathbb{N}_0$ be a non-empty set. Then, the set $S$ of functions $f \in \mathcal{A}_{F}(\Omega)$ such that all derivatives $f^{(l)}$, $l \in \mathbb{N}_0$, are unbounded on $\Omega$, is dense and $G_\delta$ in $\mathcal{A}_{F}(\Omega)$.\\\\
More generally, if $\Omega$ is an unbounded open set and $(z_n)_n$ is a sequence of points in $\Omega$ converging to $\infty$, then for $X=\{z_n:n\in \mathbb{N}\}$, generically every function $f \in \mathcal{A}_{F}(\Omega)$ has the property that the derivative $f^{(l)}$ is unbounded on $X$, for all $l \in \mathbb{N}_0$. To give an explicit example of such a function $f \in \mathcal{A}_{F}(\mathbb{C})=H(\mathbb{C})$, it suffices to set $f(z)=exp(e^{-i \theta} z)$, for some well-chosen $\theta \in \mathbb{R}$. Indeed, let $c_n=\frac{z_n}{|z_n|}$ and let $c_{k_n}{\rightarrow} c$ for a subsequence. Then, $|c|=1$ and it suffices to choose $\theta \in \mathbb{R}$, such that $c=e^{-i\theta}$. One can easily see that $|f^{(l)}(z)|=|f(z)|$, for all $z\in \mathbb{C}$ and $l \in \mathbb{N}_0$ and that $|f(z_{k_n})| \rightarrow +\infty$.\\

\noindent \textbf{Acknowledgements:} The authors would like to thank professors Ap. Giannopoulos and N. A' Campo for helpful discussions. This research was supported through the "Research in Pairs" program by the Mathematisches Forschungsinstitut Oberwolfach in 2016.\\

\section*{References and bibliography}
[1] L.V. Ahlfors, Complex Analysis: an Introduction to the Theory of Analytic Functions of One Complex Variable, 3rd Edition, McGraw-Gill, 1979\\\\\relax
[2] N. Bourbaki, Functions of a Real Variable, Springer, 1982\\\\\relax
[3] D. S. Mitrinović, Analytic Inequalities, Springer, 1970\\\\\relax
[4] G.H. Hardy, J.E. Littlewood, G.Polya, Inequalities, Cambridge University Press, 1934\\\\\relax
[5] E. Landau, Einige Ungleichungen für zweimal differenzierbare Funktionen, Proc. London Math. Soc., 13, p. 43-49, 1913\\\\\relax
[6] A.N. Kolmogorov, On inequalities between upper bounds of consecutive derivatives of an arbitrary function on an infinite interval, Amer. Math. Soc. Transl. Ser. 1, Amer. Math. Soc, Providence, p. 233–243, 1962\\\\\relax
[7] J. Garnett, Analytic Capacity and Measure, Lecture Notes in Mathematics, vol. 297, Springer, 1972\\\\\relax
[8] J. Garnett,  Positive Length but Zero Analytic Capacity, Proc. Amer. Math. Soc. 24, p. 696-699, 1970\\\\\relax
[9] S. N. Mergelyan, Certain questions of the constructive theory of functions, Trudy Mat. Inst. Steklov., 37, Acad. Sci. USSR, Moscow, p. 3–91, 1951\\\\\relax
[10] V. Nestoridis, Non-extendable holomorphic functions Math. Proc. Cambridge Phil. Sol. 139(2), p. 351-360, 2005 \\\\\relax
[11] V. Nestoridis, I. Zadik, Padé Approximants, density of rational functions in $A^\infty(\Omega)$ and smoothness of the integration operator, J.M.M.A. 423(2), p. 1534-1536, 2015, see also \href{https://arxiv.org/abs/1212.4394}{arXiv: 1212.4394}\\\\\relax
[12] W. Smith, D. M. Stolyarov, A. Volberg, Uniform approximation of Bloch functions and the boundedness of the integration operator on $H^\infty(\Omega)$, \href{https://arxiv.org/abs/1604.05433}{arXiv: 1604.05433}\\\\\relax
[13] E. Bolkas, V. Nestoridis, C. Panagiotis, M. Papadimitrakis, One sided extendability and p-continuous analytic capacities, \href{https://arxiv.org/abs/1606.05443}{arXiv: 1606.05443}\\\\\relax
[14] V. Nestoridis, Domains of Holomorphy, \href{https://arxiv.org/abs/1701.00734}{arXiv: 1701.00734}\\\\\relax
[15] M. Siskaki, Boundedness of derivatives and anti-derivatives of holomorphic functions as a rare phenomenon, \href{https://arxiv.org/abs/1611.05386}{arXiv: 1611.05386}\\\\\relax

\section*{E-mails and addresses}

Dionysios Moschonas\\
National and Kapodistrian University of Athens, Department of Mathematics\\
Panepistimiopolis, 157 84, Athens, Greece\\
e-mail: dmoschon@math.uoa.gr\\\\
Vassili Nestoridis\\
National and Kapodistrian University of Athens, Department of Mathematics\\
Panepistimiopolis, 157 84, Athens, Greece\\
e-mail: vnestor@math.uoa.gr\\\\
\end{document}